\documentclass[11pt]{amsart}
\usepackage{amsmath,amsfonts,amsthm,amsopn,amssymb}
\usepackage{cite,marginnote}

\usepackage{color,enumitem,graphicx}
\usepackage[colorlinks=true,urlcolor=blue,
citecolor=red,linkcolor=blue,linktocpage,pdfpagelabels,
bookmarksnumbered,bookmarksopen]{hyperref}
\usepackage[english]{babel}

\usepackage[left=2.75cm,right=2.75cm,top=2.9cm,bottom=2.9cm]{geometry}
\usepackage[hyperpageref]{backref}

\numberwithin{equation}{section}

\makeindex
\newtheorem{theorem}{Theorem}[section]
\newtheorem{proposition}[theorem]{Proposition}
\newtheorem{lemma}[theorem]{Lemma}
\newtheorem{remark}[theorem]{Remark}
\newtheorem{example}[theorem]{Example}
\newtheorem{corollary}[theorem]{Corollary}
\newtheorem{definition}[theorem]{Definition}
\newtheorem{theoremletter}{Theorem}

\newcommand{\ud}{\mathrm{d}}
\newcommand{\RN}{\mathbb R^N}
\newcommand{\om}{\Omega}

\newcommand{\iy}{\infty}

\newcommand{\s}{\section}
\newcommand{\dd}{\delta}
\newcommand{\DD}{\Delta}
\newcommand{\g}{\gamma}
\newcommand{\G}{\Gamma}
\newcommand{\na}{\nabla}

\newcommand{\la}{\lambda}

\newcommand{\R}{\mathbb R}
\newcommand{\al}{\alpha}

\newcommand{\ti}{\tilde}

\newcommand{\re}[1]{(\ref{#1})}

\newcommand{\rg}{\rightarrow}

\newcommand{\e}{\varepsilon}

\newcommand{\lab}{\label}
\newcommand{\bt}{\begin{theorem}}
\newcommand{\et}{\end{theorem}}
\newcommand{\bl}{\begin{lemma}}
\newcommand{\el}{\end{lemma}}
\newcommand{\bd}{\begin{definition}}
\newcommand{\ed}{\end{definition}}
\newcommand{\bc}{\begin{corollary}}
\newcommand{\ec}{\end{corollary}}
\newcommand{\bp}{\begin{proof}}
\newcommand{\ep}{\end{proof}}
\newcommand{\bx}{\begin{example}}
\newcommand{\ex}{\end{example}}
\newcommand{\bi}{\begin{exercise}}
\newcommand{\ei}{\end{exercise}}
\newcommand{\bo}{\begin{proposition}}
\newcommand{\eo}{\end{proposition}}
\newcommand{\br}{\begin{remark}}
\newcommand{\er}{\end{remark}}
\newcommand{\be}{\begin{equation}}
\newcommand{\ee}{\end{equation}}
\newcommand{\ba}{\begin{align}}
\newcommand{\ea}{\end{align}}
\newcommand{\bn}{\begin{enumerate}}
\newcommand{\en}{\end{enumerate}}
\newcommand{\bg}{\begin{align*}}
\newcommand{\bcs}{\begin{cases}}
\newcommand{\ecs}{\end{cases}}

\renewcommand{\S}{{\mathcal S}}

\newcommand{\bean}{\begin{eqnarray*}}
\newcommand{\eean}{\end{eqnarray*}}

\title[Schr\"odinger-Poisson systems at critical
growth]{Schr\"odinger-Poisson systems with \\ a general critical
nonlinearity}

\author[J. Zhang]{Jianjun Zhang}
\author[J.M.\ do \'O]{Jo\~ao Marcos do \'O}
\author[M. Squassina]{Marco Squassina}

\address[J.\ Zhang]{\newline\indent  Chern Institute of Mathematics
\newline\indent
Nankai University
\newline\indent
Tianjin 300071, PR China
\newline\indent and
\newline\indent College of Mathematics and Statistics
\newline\indent
Chongqing Jiaotong University
\newline\indent
Chongqing 400074, PR China}
\email{\href{mailto:zhangjianjun09@tsinghua.org.cn}{zhangjianjun09@tsinghua.org.cn}}

\address[J.M.\ do \'O]{\newline\indent Department of Mathematics
\newline\indent
Federal University of Para\'{\i}ba
\newline\indent
58051-900, Jo\~ao Pessoa-PB, Brazil}
\email{\href{mailto:jmbo@pq.cnpq.br}{jmbo@pq.cnpq.br}}

\address[M.\ Squassina]{\newline\indent Dipartimento di Informatica
\newline\indent
Universit\`a degli Studi di Verona,
\newline\indent
C\'a Vignal 2, Strada Le Grazie 15, I-37134 Verona, Italy}
\email{\href{mailto:marco.squassina@univr.it}{marco.squassina@univr.it}}

\thanks{Research partially supported by INCTmat/MCT/Brazil.\
J.M.\ do \'O was supported by CNPq, CAPES/Brazil,
Jianjun Zhang  was partially supported by CAPES/Brazil and CPSF (2013M530868)}

\subjclass[2000]{35B25, 35B33, 35J61}
\keywords{Schr\"odinger-Poisson systems, variational methods, critical growth}

\begin{document}

\begin{abstract}
We consider a Schr\"{o}dinger-Poisson system involving a
general nonlinearity at critical growth and
we prove the existence of positive solutions. The
Ambrosetti-Rabinowitz condition is not required. We
also study the asymptotics of solutions with respect to a parameter.
\end{abstract}
\maketitle

\s{Introduction and main result}

\noindent
We are concerned with the nonlinear
Schr\"{o}dinger-Poisson system
\be\lab{q1} \left\{
\begin{array}{ll}
-\Delta u+u+\la\phi u=f(u)\ \ \ \mbox{in}\ \R^3,\\
-\DD\phi=\la u^2,\ \ \ \ \ \ \ \ \ \ \ \ \ \ \ \ \ \mbox{in}\ \R^3,
\end{array}
\right. \ee where $\la>0$ and the nonlinearity $f$ reaches the critical
growth. In the last decade, the Schr\"{o}dinger-Poisson system has
been object of intensive research because of its strong
relevance in applications. From a physical point of view, it describes systems of
identically charged particles interacting each other in the case where
magnetic effects can be neglected.
The nonlinear term $f$ models the
interaction between the particles and the coupled term $\phi u$
concerns the interaction with the electric field. For more detailed
physical aspects of the Schr\"{o}dinger-Poisson system, we refer the reader to
\cite{Ruiz1,Benci1,Benci2,Vaira} and to the references therein.
In recent years, there has been an increasing attention towards
systems like \re{q1} and the existence of positive solutions,
sign-changing solutions, ground states, radial and non-radial
solutions and semi-classical states has been investigated.
In \cite{T1}, D'Aprile and Mugnai obtained the existence of a nontrivial radial solution to
\re{q1} with $f(u)=|u|^{p-2}u$, for $p\in[4,6)$. In \cite{T2},
D'Aprile proved that system \re{q1} admits a non-radial solution
for $f(u)=|u|^{p-2}u$, with $p\in(4,6)$. In \cite{AP}, by using the
Concentration Compactness Principle, Azzollini and Pomponio obtained
the existence of a ground state solution to \re{q1} with
$f(u)=|u|^{p-2}u$, for $p\in(3,6)$. In \cite{Ruiz}, Ruiz obtained some
nonexistence results for
\be\lab{qqqq1} \left\{
\begin{array}{ll}
-\Delta u+u+\lambda\phi u=|u|^{p-2}u\ \ \ \mbox{in}\ \R^3,\\
-\DD\phi=u^2,\ \ \ \ \ \ \ \ \ \ \ \ \ \ \ \ \ \ \ \ \ \ \
 \mbox{in}\ \R^3
\end{array}
\right. \ee and established the relation between the existence of the positive solutions to system \re{qqqq1} and the parameters
$p\in(2,6)$ and $\lambda>0$. Moreover, if $\la\ge\frac{1}{4}$, the author
showed that $p=3$ is a critical value for the existence of the
positive solutions. For $p\in(2,3)$, Ruiz \cite{R2} investigated the
existence of radial ground states to system \re{qqqq1} and
obtained the different behavior of the solutions depending on $p$ as
$\la\rg 0$. We also would like to cite some works \cite{T3,R1},
where system \re{qqqq1} was considered as $\la\rg 0$. In
\cite{T3,R1}, the authors were concerned with the semi-classical
states for system \re{qqqq1}. Precisely, the authors studied the
existence of radial positive solutions concentrating around a
sphere. Recently, some works were focused on the existence of
sign-changing solutions to \re{q1} with $f(u)=|u|^{p-2}u$.
By using a gluing method, Kim and Seok \cite{Kim} proved the
existence of sign-changing solutions with a prescribed number of
nodal domains for \re{q1} with $p\in(4,6)$. Subsequently,
Ianni \cite{Ianni1} obtained a similar result for $p\in[4,6)$. More
recently, Wang and Zhou \cite{WZ} considered the
non-autonomous system
\be\lab{qqqq2} \left\{
\begin{array}{ll}
-\Delta u+V(x)u+\lambda\phi u=|u|^{p-2}u\ \ \ \mbox{in}\ \R^3,\\
-\DD\phi=u^2,\ \ \ \ \ \ \ \ \ \ \ \ \ \ \ \ \ \ \ \ \ \ \ \ \ \ \ \
\
 \mbox{in}\ \R^3.
\end{array}
\right. \ee Under suitable conditions on $V$, they proved
the existence of least energy sign-changing solutions to system
\re{qqqq2} with $p\in(4,6)$ by minimizing over the sign-changing
Nehari manifold. For further works on the non-autonomous
Schr\"{o}dinger-Poisson system, we also would like to mention
\cite{Am,CV,JZ,LP,ZZ} and the references therein.

The works discussed above mainly focus on the study of system
\re{q1} with the very special nonlinearity $f(u)=|u|^{p-2}u$. In
\cite{Azzo}, Azzollini, d'Avenia and Pomponio were concerned with
the existence of a positive radial solution to system \re{q1}
under the effect of a general nonlinear term, see also
\cite{AzzdAvpomp,jeongseok}. Precisely, let $g(u)=-u+f(u)$, then

\begin{theoremletter}[see
\cite{Azzo}] {\it Suppose
\begin{itemize}
\item [$(H1)$] $g(s)\in C(\mathbb{R},\mathbb{R})$;
\item [$(H2)$] $\displaystyle -\infty<\liminf_{s\rightarrow0}\frac{g(s)}{s}\le
\limsup_{s\rightarrow0}\frac{g(s)}{s}=-m<0$;
\item [$(H3)$] $\displaystyle
\limsup_{s\rightarrow\infty}\frac{g(s)}{s^5}\le0$;
\item [$(H4)$] there exists $\xi_0>0$ such that $G(\xi_0):=\int_0^\xi
g(s)\, \ud s>0$.
\end{itemize} Then there exists $\la_0>0$ such that
\re{q1} admits a positive radial solution for
$\la\in(0,\la_0)$.}
\end{theoremletter}

\noindent
$(H_1)$-$(H_4)$ are known as Berestycki-Lions
conditions, introduced in \cite{Lions}. There, the authors showed that these
conditions are almost necessary and sufficient for
the existence of ground states to the nonlinear scalar field
equation $-\DD u=g(u)$, with $u\in H^1(\RN)$, $N\ge 3$.

\noindent
We remark that in the literature described above, only the {\em subcritical} case
was considered. A natural question arises on  whether results
like Theorem A holds if $f$ is at {\em critical growth}. In fact, in \cite{Jian},
Zhang, obtained the following


\begin{theoremletter}[see \cite{Jian}] {\it Suppose $f\in C(\mathbb{R},\mathbb{R})$ is odd and
\begin{itemize}
\item [($g_1$)] $\lim\limits_{s\rightarrow 0}\frac{f(s)}{s}=0$,
\item [($g_2$)] $\lim\limits_{s\rightarrow +\infty}\frac{f(s)}{s^5}=K>0$,
\item [($g_3)$] There exists $D>0$ and $q\in(2,6)$ such that $f(s)\ge Ks^5+D s^{q-1}$, for all $s>0$,
\item [($g_4)$] There exists $\gamma>2$ such that $0<\gamma \int_0^sf(\tau)\, \ud \tau \le sf(s)$, for all $s\not=0$,
\end{itemize}
Then $(i)$ \re{q1} has a positive radial solution for small $\la>0$ if $q\in(2,4]$ with $D$ large enough, or $q\in(4,6)$;
$(ii)$ if $\gamma>3$, \re{q1} admits a ground state solution for any $\la>0$ provided $q\in(2,4]$ with $D$ large enough, or $q\in(4,6)$.
}
\end{theoremletter}

\noindent
The author was able to obtain this existence result
via a truncation argument. Condition $(g_4)$ implies that $f$ is superlinear and is the
so-called Ambrosetti-Rabinowitz condition, usually involved in guaranteeing the
boundedness of (PS)-sequences.

\vskip3pt
\noindent
The aim of this paper is to study the existence of the positive
solutions to system \re{q1} involving a more general critical
nonlinearity compared to that allowed in Theorem B.
In particular, the Ambrosetti-Rabinowitz condition is not required.
\vskip2pt
\noindent
We shall assume that the following hypotheses on $f$:
\begin{itemize}
\item [($f_1$)] $f:\R\rg \R$ is continuous, $f=0$ on $\R^-$ and $\lim\limits_{s\rightarrow 0}\frac{f(s)}{s}=0$.
\item [($f_2$)] $\limsup\limits_{s\rightarrow +\infty}\frac{f(s)}{s^5}\le 1$.
\item [($f_3)$] There exist $\mu>0$ and $q\in(2,6)$ such that $f(s)\ge \mu s^{q-1}$ for all $s\geq 0$.
\end{itemize}

\noindent
Assumption $(f_2)$ implies $f$ has (possibly) a critical growth at
infinity and the limit of $f(s)/s^5$ at $+\infty$
may fail to exist. Moreover, there exists $\kappa>0$ such that
\be
\lab{f}
f(s)\le \frac{1}{2}s+\kappa s^5\ \ \ \mbox{for all $s\ge 0$}.
\ee
Before stating the main result, we fix some notations. In
the sequel, $\S$ and $\mathcal{C}_q$ denote the best constants of
Sobolev embeddings $\mathcal{D}^{1,2}(\R^3)\hookrightarrow L^{6}(\R^3)$ and
$H^1(\R^3)\hookrightarrow L^q(\R^3),$
\begin{align*}
 \S\left(\int_{\R^3}|u|^6\, \ud x \right)^{\frac{1}{3}} &\le \int_{\R^3}|\na
u|^2\, \ud x ,\quad \text{for all $u\in \mathcal{D}^{1,2}(\R^3)$}, \\
 \mathcal{C}_q\left(\int_{\R^3}|u|^q\, \ud x \right)^{\frac{2}{q}} &\le
\int_{\R^3}(|\na u|^2+|u|^2)\, \ud x ,\quad \text{for all $u\in H^1(\R^3)$}.
\end{align*}
Our main result is the following
\bt\lab{Th1} Suppose that $f$ satisfies $(f_1)$-$(f_3)$.
\begin{itemize}
\item [(i)] There exists $\la_0>0$ such that, for every $\la\in (0,\la_0)$, system
\re{q1} admits a nontrivial positive solution
$(u_\la,\phi_\la)$, provided that
$$
\mu>\left[\frac{3q-6}{2q\S^\frac{3}{2}}\right]^{\frac{q-2}{2}}\mathcal{C}_q^{\frac{q}{2}};
$$
\item [(ii)] Along a subsequence, $(u_\la,\phi_\la)$ converges to $(u,0)$ in
$H^1(\R^3)\times\mathcal{D}^{1,2}(\R^3)$ as $\la\rg 0$, where $u$ is
a ground state solution to the limit problem
$$
-\DD u+u=f(u),\ \ \ u\in H^1(\R^3).
$$
\end{itemize}
\et
\noindent
The rest of the paper is devoted to prove Theorem \ref{Th1}. Since
we are concerned with system \re{q1} with a more general
nonlinear term $f$, the problem becomes more thorny and tough in
applying variational methods. In fact, due to the lack of the
Ambrosetti-Rabinowitz condition, the boundedness of (PS)-sequence is
not easy to be obtained. To overcome this difficulty, we will adopt
a local deformation argument from Byeon and Jeanjean \cite{byeon} to
get a bounded (PS)-sequence. Due to the presence of the
nonlocal term $\phi u$, a crucial modification on the min-max
value is needed. We will define another min-max value $C_\la$ (see
Section 3), where all paths are required to be uniformly bounded
with respect to $\la$. Similar arguments can be found in \cite{Chen-Zou}.

\vskip0.1in

The paper is organized as follows. \newline
In Section 2 we consider the functional
framework and some preliminary results. \newline
In Section 3 we construct the min-max level. \newline
In Section 4, we use a local deformation argument to give the proof of Theorem~\ref{Th1}.

\vskip0.12in

\noindent{\bf Notations.}
\begin{itemize}
\item [$\bullet$] $\|u\|_p:=\big(\int_{\R^3}|u|^p\, \ud x\big)^{1/p}$ for $p\in [1,\infty)$.
\item [$\bullet$] $\|u\|:=\big(\|u\|^2_2+\|\nabla u\|^2_2\big)^{1/2}$ for $u\in H^1(\R^3)$.
\item [$\bullet$] $H_r^1(\R^3)$ is the subspace of $H^1(\R^3)$ of radially symmetric functions.
\item [$\bullet$] $\mathcal{D}^{1,2}(\R^3): =\{u\in L^{2^{\ast}}(\R^3): \nabla u\in L^2(\R^3)\}$.
\end{itemize}


\s{Preliminaries and functional setting}

\renewcommand{\theequation}{2.\arabic{equation}}

\noindent
We recall that, for $u\in H^1(\R^3)$, the Lax-Milgram theorem
implies that there exists a unique $\phi_u\in
\mathcal{D}^{1,2}(\R^3)$ such that $-\DD\phi=\la u^2$ with
\be
\lab{t1} \phi_u(x):=\la\int_{\R^3}\frac{u^2(y)}{4\pi|x-y|}\, \ud y.
\ee
Setting
$$
T(u):=\frac{1}{4}\int_{\R^3}\phi_uu^2\, \ud x,
$$
we summarize some properties of $\phi_u, T(u)$,
which will be used later.

\bl[see \cite{Ruiz}]
\lab{l0} For any $u\in H^1(\R^3)$, we have
\begin{itemize}
\item [(1)] $\phi_u:H^1(\R^3)\mapsto \mathcal{D}^{1,2}(\R^3)$ is continuous and maps bounded sets into bounded
sets.
\item [(2)] $\phi_u\ge0, T(u)\le
c\la\|u\|^4$ for some $c>0$.
\item [(3)] If $u_n\rg u$ weakly in $H^1(\R^3)$, then $\phi_{u_n}\rg \phi_u$ weakly in
$\mathcal{D}^{1,2}(\R^3)$.
\item [(4)] If $u_n\rg u$ weakly in $H^1(\R^3)$, then
$T(u_n)=T(u)+T(u_n-u)+o(1)$.
\item [(5)] If $u$ is a radial function, so is $\phi_u$.
\end{itemize}
\el
\noindent
Substituting \re{t1} into \re{q1}, we can rewrite
\re{q1} in the following equivalent equation
\be\lab{q2} -\DD u+u+\la\phi_uu=f(u),\ \ \ u\in H^1(\R^3). \ee We
define the energy functional $\G_\la:H^1(\R^3)\to\R$ by
$$
\G_\la(u)=\frac{1}{2}\int_{\R^3}(|\na
u|^2+u^2)\, \ud x+\frac{\la}{4}\int_{\R^3}\phi_uu^2\, \ud x-\int_{\R^3}F(u)\, \ud x,
$$
with $F(t)=\int_0^tf(s)\, \ud s$. It is standard to show that $\G_\la$ is of class $C^1$ on $H^1(\R^3)$. Since we are concerned with the positive solutions
of \re{q1}, from now on, we can assume that $f(s)=0$ for every $s\leq 0$.
It is readily proved that any critical point of $\G_\la$ is
nonnegative and, by the maximum principle, it is strictly positive.
Moreover, it is easy to verify that $(u,\phi)\in
H^1(\R^3)\times\mathcal{D}^{1,2}(\R^3)$ is a solution of \re{q1}
if and only if $u\in H^1(\R^3)$ is a critical point of the
functional $\G_\la$. If $\la=0$, problem \re{q2} becomes
\be
\lab{q3} -\DD u+u=f(u),\ \ \ u\in H^1(\R^3),
\ee
which will be referred as
the limit problem of \re{q2}.\ In general, if a problem is well-behaved and undergoes a
small perturbation, then one may expect that the perturbed problem
has a solution near the solutions of the original problem. Then if
$\la$ is small, it is natural to find a solution of \re{q2} in
some neighborhood of the solutions to the limit problem \re{q3}, which will play a crucial r\v ole in the
study of perturbed problem \re{q2}.
In the following, we study some properties of
the limit problem \re{q3}. First, we show the existence of the
ground states of the limit problem \re{q3}.

\bo\lab{bo1} Suppose that $f$ satisfies $(f_1)$-$(f_3)$, then the
limit problem \re{q3} has a ground state $u\in H_r^1(\R^3)$,
provided that
\begin{equation}
 \label{mucond}
 \mu>\left[\frac{3q-6}{2q\S^\frac{3}{2}}\right]^{\frac{q-2}{2}}\mathcal{C}_q^{\frac{q}{2}}.
 \end{equation}
\eo
\br\rm  In \cite{souto} the authors established the existence of the
ground state solutions for the nonlinear scalar field equation
involving critical growth in $\RN$, for $N\ge 2$. In particular, assuming
that $f$ satisfies $(f_1)$-$(f_3)$ and an additional condition
\begin{itemize}
\item [($f_4)$] $sf(s)-2F(s)\ge 0$ for all
$s\ge0$, where $F(s)=\int_0^sf(\tau)\, \ud \tau$,
\end{itemize}
the authors proved that \re{q3} has a ground state. We
remark that $(f_4)$ can be removed. \er

\vskip3pt
\noindent
To prove Proposition \ref{bo1}, we will use the following notations.
\begin{align*}
\mathcal{M}&:=\left\{u\in H_r^1(\R^3)\setminus \{0\}:
\int_{\R^3}G(u)\, \ud x=1\right\}, \\
\mathcal{P}&:=\left\{u\in H_r^1(\R^3)\setminus \{0\}:
6\int_{\R^3}G(u)\, \ud x=\int_{\R^3}|\na u|^2\, \ud x\right\},
\end{align*}
where $G(t)=F(t)-\frac{1}{2}t^2$. $\mathcal{P}$ is the so-called
Pohoz$\check{a}$ev manifold. It follows, from $(f_3)$, that there
exists $\xi>0$ such that $G(\xi)>0$. Then it is easy to check that
$\mathcal{M}\neq\emptyset$ and $\mathcal{P}\neq\emptyset$. Define
$$
M:=\frac{1}{2}\inf_{u\in\mathcal{M}}\int_{\R^3}|\na u|^2\, \ud x,\quad\,\,
p:=\inf_{u\in\mathcal{P}}I(u),
$$
and the Mountain Pass value
$$
b:=\inf_{\g\in\Upsilon}\max_{0\le t\le 1}I(\g(t)),
$$
where $\Upsilon=\{\g\in C([0,1],H_r^1(\R^3)):\g(0)=0,\, I(\g(1))<0\}$
and
$$
I(u):=\frac{1}{2}\int_{\R^3}(|\na u|^2+|u|^2)\, \ud x-\int_{\R^3}F(u)\, \ud x.
$$
\bl\lab{l1}
Let $f$ satisfy $(f_1)$-$(f_3)$ and \eqref{mucond}.
Then $0<M<\frac{\sqrt[3]{6}}{2}\S$ and $p<\frac{1}{3}\S^{\frac{3}{2}}$.
\el
\bp\noindent
Obviously $M\in[0,\iy)$. We claim that $M>0$. Assume by contradiction that it is $M=0$. Then
there exists $\{u_n\}_n\subset\mathcal{M}$ such that $\|\na u_n\|_2\rg 0$ as $n\rg\iy$.
By Sobolev's embedding theorem, $\|u_n\|_6\rg0$ as $n\rg\iy$. Thus, it follows from \re{f} that
$$
\limsup_{n\rg\iy}\int_{\R^3}G(u_n)\, \ud x\le\limsup_{n\rg\iy}\frac{\kappa}{6}\int_{\R^3}|u_n|^6\, \ud x=0,
$$
a contradiction, proving the claim.
We now claim that $p\le b$. It suffices to prove that
$$
\g([0,1])\cap
\mathcal{P}\neq \emptyset\ \ \mbox{for all $\g\in \Upsilon$},
$$
whose proof is similar to that in \cite[Lemma 4.1]{Jean}. Let
$$
P(u)=\int_{\R^3}|\na u|^2\, \ud x-6\int_{\R^3}G(u)\, \ud x.
$$
Then by \re{f} it is easy to know that there exists $\rho_0>0$ such that
\be\lab{ss12}
P(u)>0\ \ \mbox{if}\ \ 0<\|u\|\le \rho_0.
\ee
For any $\g\in\Upsilon$, $P(\g(0))=0$ and $P(\g(1))\le 6 I(\g(1))<0$. Thus, there exists $t_0\in(0,1)$ such that $P(\g(t_0))=0$ with $\|\g(t_0)\|>\rho_0$, which implies $\g([0,1])\cap\mathcal{P}\neq \emptyset$.
We now use an idea from Coleman-Glazer-Martin
\cite{Coleman} to prove that $p=\frac{2\sqrt{3}}{9}M^{\frac{3}{2}}$.
Define $\Phi: \mathcal{M}\rightarrow
\mathcal{P}$: $(\Phi(u))(x)=u(\frac{x}{t_u})$, where
$t_u=\sqrt{6}/6\|\nabla u\|_2$. Then $\Phi$ is a bijection. For $u\in H^1(\R^3)$, let us set
$$
T_0(u)=\frac{1}{2}\int_{\R^3}|\nabla u|^2\, \ud x,
\,\,\quad V(u)=\int_{\R^3}G(u)\, \ud x.
$$
Then for
$u\in \mathcal{M}$,
$I(\Phi(u))=t_uT_0(u)-t_u^3V(u)=\sqrt{6}/18\|\nabla
u\|_2^3$. Thus, $$\inf_{u\in \mathcal{P}}I(u)=\inf_{u\in
\mathcal{M}}I(\Phi(u))=
\sqrt{6}/18\inf_{u\in\mathcal{M}}\|\nabla
u\|_2^3,$$
which implies $p=\frac{2\sqrt{3}}{9}M^{\frac{3}{2}}$.
Finally, similar to that in \cite{souto}, taking $\psi\in H_r^1(\R^3)$ with $\psi\ge 0$ with $\|\psi\|_q^2=\mathcal{C}_q^{-1}$ and $\|\psi\|=1$,
then
\begin{equation*}
b\le\max_{t\ge0}I(t\psi)\le\max_{t\ge0}
\Big(\frac{t^2}{2}-\mu\frac{t^q}{q}\|\psi\|_q^q\Big)
=\frac{q-2}{2q}\mu^{-\frac{2}{q-2}}\mathcal{C}_q^{\frac{q}{q-2}}.
\end{equation*}
Thus, by virtue of \eqref{mucond}, we have $p<\frac{1}{3}\S^{\frac{3}{2}}$ and, in turn,
$M<\frac{\sqrt[3]{6}}{2}\S$.
\ep

\noindent
In the following, we will show that $p$ can be achieved. This implies
that the limit problem \re{q3} admits a ground state solution.
Similar to that in \cite{Lions}, it is enough to prove that $M$ can
be achieved. Now, we give the following Brezis-Lieb Lemma.

\bl\lab{l2} Let $h\in C(\R^3\times\R)$ and suppose that
\be\lab{rr}\lim_{t\rg 0}\frac{h(x,t)}{t}=0\ \ \mbox{and}\ \
\limsup_{|t|\rg \iy}\frac{|h(x,t)|}{|t|^5}<\iy, \ee uniformly in
$x\in\R^3$. If $u_n\rg u_0$ weakly in $H^1(\R^3)$ and $u_n\rg u_0$
$a.e.$ in $\R^3$, then
$$
\int_{\R^3}(H(x,u_n)-H(x,u_n-u_0)-H(x,u_0))\, \ud x=o(1),
$$
where $H(x,t)=\int_0^th(x,s)\, \ud s$. \el
\begin{proof}
The proof is standard. For the subcritical case, we refer to V. Coti Zelati and P.H. Rabinowitz \cite{Coti}. For any fixed $\dd>0$, set $\om_n(\dd):=\{x\in\R^3: |u_n(x)-u_0(x)|\le \dd\}$. Then
\begin{align*}
&\int_{\R^3}(H(x,u_n)-H(x,u_n-u_0)-H(x,u_0))\, \ud x\\
=&\int_{\R^3\setminus\om_n(\dd)}(H(x,u_n)-H(x,u_n-u_0)-H(x,u_0))\, \ud x\\
&+\int_{\om_n(\dd)}(H(x,u_n)-H(x,u_0))\, \ud x-\int_{\om_n(\dd)}H(x,u_n-u_0)\, \ud x\\
:=&J_1+J_2+J_3.
\end{align*}
By conditions \re{rr}, for any $\rho>0$, there exists $C_\rho>0$ such that
$|h(x,t)|\le \rho |t|+C_\rho|t|^5$ for all $(x,t)\in \R\times\R^3$.
Then
\begin{align*}
|J_3|&\le\int_{\om_n(\dd)}\left(\frac{\rho}{2}|u_n-u_0|^2+\frac{C_\rho}{6}|u_n-u_0|^6\right)\, \ud x\\
&\le\left(\frac{\rho}{2}+\frac{C_\rho}{6}\dd^4\right)
\int_{\R^3}|u_n-u_0|^2\, \ud x,
\end{align*}
and {\allowdisplaybreaks
\begin{align*}
|J_2|&\le\int_{\om_n(\dd)}\left[\rho(|u_n|+|u_0|)+C_\rho(|u_n|+|u_0|)^5\right]|u_n-u_0|\, \ud x\\
&\le\rho\left(\int_{\R^3}(|u_n|+|u_0|)^2\, \ud x\right)^{\frac{1}{2}}
\left(\int_{\R^3}|u_n-u_0|^2\, \ud x\right)^{\frac{1}{2}}\\
&+C_\rho\left(\int_{\R^3}(|u_n|+|u_0|)^6\, \ud x\right)^{\frac{5}{6}}
\left(\int_{\om_n(\dd)}|u_n-u_0|^6\, \ud x\right)^{\frac{1}{6}}\\
\le&\rho\left(\int_{\R^3}(|u_n|+|u_0|)^2\, \ud x\right)^{\frac{1}{2}}
\left(\int_{\R^3}|u_n-u_0|^2\, \ud x\right)^{\frac{1}{2}}\\
&+C_\rho\dd^{\frac{2}{3}}\left(\int_{\R^3}(|u_n|+|u_0|)^6\, \ud x\right)^{\frac{5}{6}}
\left(\int_{\R^3}|u_n-u_0|^2\, \ud x\right)^{\frac{1}{6}}.
\end{align*}
}%
Then, since $\{u_n\}_n$ is bounded in $H^1(\R^3)$, for every $\e>0$, there
exist $\rho, \dd>0$ such that $|J_2|+|J_3|\le \e/2$, for all $n\geq 1$.
On the other hand,
\begin{align*}
J_1&=\int_{B_R(0)\setminus\om_n(\dd)}(H(x,u_n)-H(x,u_n-u_0)-H(x,u_0))\, \ud x\\
&\ \ \ \ +\int_{\R^3\setminus(\om_n(\dd)\bigcup
B_R(0))}(H(x,u_n)-H(x,u_n-u_0)-H(x,u_0))\, \ud x\\
&:=K_1+K_2,
\end{align*}
where $B_R(0)=\{x\in\RN:|x|<R\}, R>0$. Noting that
\begin{align*}
|K_2|\le\int_{\R^3\setminus
B_R(0)}\left[\rho(|u_n|+|u_0|)+C_\rho(|u_n|+|u_0|)^5\right]|u_0|\, \ud x
+\int_{\R^3\setminus B_R(0)}H(x,u_0)\, \ud x,
\end{align*}
there exists $R>0$ with $|K_2|\le \e/4$, for all $n\geq 1$.
Recall that $u_n\rg u_0$ a.e.\ in $\R^3$, then it follows from the Severini-
Egoroff theorem that $u_n$ converges to $u_0$ in measure in $B_R(0)$, which imples
$$
\lim\limits_{n\rg 0}|B_R(0)\setminus\om_n(\dd)|=0.
$$
In turn $|K_1|\le \e/4$ for $n$ large. Then $|J_1|\le \e/2$ for $n$ large and
the proof is complete.
\end{proof}
\vskip2pt

\noindent\textit{\bf Proof of Proposition \ref{bo1}.}
The proof is similar to that of \cite{Zhang-Zou}. We may assume that there exists
$\{u_n\}\subset H_r^1(\R^3)$ such that $\int_{\R^3}G(u_n)\, \ud x=1$ and
$\int_{\R^3}|\na u_n|^2\, \ud x\rg 2M,$ as $n\rg\iy$. By $(f_1)$-$(f_2)$,
$\{u_n\}$ is bounded in $H_r^1(\R^3)$. Thus there is $u_0\in H_r^1(\R^3)$ such that,
up to a subsequence, $u_n\rg u_0$ weakly in $H_r^1(\R^3)$. Then
$$
\int_{\R^3}|\na u_n|^2\, \ud x=\int_{\R^3}|\na u_0|^2\, \ud x+\int_{\R^3}|\na
u_n-\na u_0|^2\, \ud x+o(1).
$$
Moreover, by Lemma \ref{l2}, we have
$$
\int_{\R^3} G(u_n)\, \ud x=\int_{\R^3} G(u_0)\, \ud x+\int_{\R^3} G(u_n-u_0)\, \ud x+o(1).
$$
It is easy to know that $M=\inf\{T_0(u):V(u)=1,u\in H_r^1(\R^3)\}$. Moreover,
$$
T_0(u_n)=T_0(v_n)+T_0(u_0)+o(1),
\quad\,
V(u_n)=V(v_n)+V(u_0)+o(1),
$$
where $v_n=u_n-u_0$. Set $S_n=T_0(v_n), S_0=T_0(u_0), V(v_n)=\la_n, V(u_0)=\la_0$, we have $\la_n=1-\la_0+o(1)$ and $S_n=M-S_0+o(1)$.
To prove that $u_0$ is a minimizer of $M$, it suffices to prove  $\lambda_0=1$, which implies $u_n\rg u_0$ strongly in $H_r^1(\R^3)$. It is easy to see that
\be\lab{qqq4}
T_0(u)\ge M(V(u))^{1/3},
\ee
for all $u\in H^1(\R^3)$ and $V(u)\ge0$. As we can see in \cite{Zhang-Zou}, $\lambda_0\in [0,1]$. If $\lambda_0\in [0,1)$, then $\lambda_n>0$ for $n$ large enough. By \re{qqq4}, we have that $S_0\ge M (\lambda_0)^{1/3}$ and  $S_n\ge M (\lambda_n)^{1/3}$. This implies
\begin{align*}
M&=\lim_{n\rightarrow \infty}(S_0+S_n)\ge\lim_{n\rightarrow \infty}M\left((\lambda_0)^{1/3}+(\lambda_n)^{1/3}\right)\\
&=M\left((\lambda_0)^{1/3}+(1-\lambda_0)^{1/3}\right)
\ge M(\lambda_0+1-\lambda_0)=M,
\end{align*}
which implies that $\lambda_0=0$. So we get that $u_0=0$ and $\lim_{n\rightarrow \infty}S_n=M$. By $(f_1)$-$(f_2)$, for any $\e>0$, there exists $C_\e>0$ such that $F(s)\le \frac{1}{4}s^2+C_\e s^4+(1+\e)s^6/6$ for $s\in\R$. Then
$$
1=\lim_{n\rg\infty}\lambda_n\le\frac{1+\e}{6}\limsup_{n\rg\iy}\|v_n\|_6^6,
$$
since $\|v_n\|_4\to 0$, namely
$\limsup_{n\rightarrow\infty}\|v_n\|_6^2\ge \sqrt[3]{6}$. Thus
$$
M=\frac{1}{2}\limsup_{n\rightarrow \infty}\|\nabla v_n\|_2^2\ge \frac{\S}{2}\limsup_{n\rightarrow\infty}\|v_n\|_6^2
\ge\frac{\sqrt[3]{6}}{2}\S,
$$
which contradicts Lemma \ref{l1}. Therefore, we conclude that $\lambda_0=1$. Therefore, $u_0\in \mathcal{M}$ and
$\int_{\R^3}|\na u_0|^2\, \ud x=2M$. Setting $t_0=\|\nabla
u_0\|_2/\sqrt{6}$, it follows from Coleman-Glazer-Martin
\cite{Coleman} that $\omega=u_0(\frac{\cdot}{t_0})\in \mathcal{P}$
is a ground state solution to problem \re{q3}. \qed

\vskip3pt
\noindent
Define as $\mathcal{S}_r$ the set of the {\em radial ground states}
$U$ of \re{q3}. Then $\omega\in\mathcal{S}_r$. Moreover, thanks
to Lemma \ref{l2}, similarly as that in \cite{Byeon-Zhang-Zou,Jean,Zhang-Zou}, we
have the following

\bo\lab{bo2}\noindent
\begin{itemize}
\item [$(i)$] $b=I(\omega)$, namely the Mountain Pass value agrees with the least energy level.
\item [$(ii)$] $\mathcal{S}_r$ is compact in $H_r^1(\R^3)$.
\end{itemize}
\eo
\begin{proof}
$(i)$ Obviously, by $(f_1)$-$(f_3)$ we know that $b$ is well defined. As we can see in the proof of Lemma \ref{l1} and Proposition \ref{bo1}, we get that $p\le b$ and $p=I(\omega)$. To prove $b$ is the least energy, it suffices to prove $b\le I(\omega)$. Noting that $\omega$ is a ground state solution to \re{q3}, similar to that in \cite{Jean}, there exists a path
$\gamma\in \Gamma$ satisfying $\gamma(0)=0, I(\gamma(1))<0,
\omega\in \gamma([0,1])$ and $\max_{t\in[0,1]}I(\gamma(t))=I(\omega)$. Thus, $b\le I(\omega)$.

$(ii)$ We adopt some ideas in \cite{Byeon-Zhang-Zou} to show the compactness of $\mathcal{S}_r$. Similar to that in \cite{Byeon-Zhang-Zou}, $\mathcal{S}_r$ is bounded in $H_r^1(\R^3)$. For any $\{u_n\}\subset\mathcal{S}_r$, without loss of generality, we can assume that $u_n\rg u_0$ weakly in $H_r^1(\R^3)$ and $u_n\rg u_0$ a.e. in $\R^3$. It follows from \cite{Maris} that $v_n(\cdot):=u_n(\sqrt{M/3}\, \cdot)$ is a minimizer of $T_0(v)$ on $\{v\in H_r^1(\R^3): V(v)=1\}$. This means that $\{v_n\}$ is a positive and radially symmetric minimizing sequence of $M$. As we can see in the proof of Proposition \ref{bo1}, $v_n\rg v_0:=u_0(\sqrt{M/3}\, \cdot)$ strongly in $H_r^1(\R^3)$. Thus, $u_n\rg u_0$ strongly in $H_r^1(\R^3)$ and $u_0\in \mathcal{S}_r$, i.e., $\mathcal{S}_r$ is compact.
\end{proof}

\vskip0.1in


\s{The minimax level}
\renewcommand{\theequation}{3.\arabic{equation}}

\noindent
Let $U\in\mathcal{S}_r$ be arbitrary but fixed. By the Pohoz\v{a}ev
identity, for $U_t(x)=U(\frac{x}{t})$ we have
$$
I(U_t)=\Big(\frac{t}{2}-\frac{t^3}{6}\Big)\int_{\R^3}|\nabla
U|^2\, \ud x.
$$
Thus, there exists $t_0>1$ such that $I(U_t)<-2$ for $t\ge t_0$. Set
$$
D_\la\equiv\max_{t\in[0,t_0]}\G_\la(U_t).
$$
Then, by virtue of Lemma \ref{l0}, we get that $D_\la\rg b$, as $\la\rg 0$.

\noindent
Moreover, it is easy to verify the following lemma, which is crucial
to define the uniformly bounded set of the mountain pathes as
previously mentioned.
\bl\lab{l4.1} There exist $\la_1>0$ and  $\mathcal{C}_0>0$, such that for
any $0<\la<\la_1$ there hold
\begin{equation*}
\G_\la(U_{t_0})<-2,\qquad
\|U_t\|\le \mathcal{C}_0, \,\,\, \forall t\in (0,t_0],\qquad
\|u\|\le \mathcal{C}_0, \,\,\, \forall u\in\mathcal{S}_r.
\end{equation*}
\el
\bp
Due to the Pohoz$\check{a}$ev identity, as we can see in \cite{Byeon-Zhang-Zou}, there exists $C>0$ such that $\|u\|\le C$ for any $u\in\mathcal{S}_r$. For $U\in\mathcal{S}_r$ fixed above and $t\in (0,t_0]$,
$$
\|U_t\|^2=t\|\na U\|_2^2+t^3\|U\|_2^2\le(t+t^3)\|U\|^2\le C^2(t_0+t_0^3).
$$
The second and last part of the assertion hold if  $\mathcal{C}_0=2t_0^2C$. For the first part, by Lemma \ref{l0}
$$
\G_\la(U_{t_0})\le I(U_{t_0})+4c\la^2\|U_{t_0}\|^4\le I(U_{t_0})+4c\la^2\mathcal{C}_0^4.
$$
It follows from $I(U_{t_0})<-2$ that there exist $\la_1>0$ with $\G_\la(U_{t_0})<-2$ for
any $0<\la<\la_1$. The proof is completed.
\ep
\noindent
Now, for any $\la\in (0,\la_1)$, we define a min-max value $C_\la$:
$$C_\la=\inf_{\g\in \Upsilon_\la}\max_{s\in [0,t_0]}\G_\la(\g(s)),$$
where
\begin{align*}
\Upsilon_\la=\big\{\g\in C([0,t_0],H_r^1(\R^3)): \g(0)=0,
\g(t_0)=U_{t_0},\|\g(t)\|\le \mathcal{C}_0+1,t\in[0,t_0]\big\}.
\end{align*}
Obviously, $U_t\in \Upsilon_\la$. Moreover, $C_\la\le D_\la$ for
$\la\in (0,\la_1)$.

\bo\lab{bo4} $\lim\limits_{\la\rg 0}C_\la=b.$ \eo

\bp It suffices to prove that
$$
\liminf\limits_{\la\rg 0}C_\la\ge b.
$$
Noting that $\phi_u\ge 0$, we see that for any $\g\in \Upsilon_\la$,
$\ti{\g}(\cdot)=\g(t_0\cdot)\in \Upsilon$. It follows that $C_\la\ge
b$, concluding the proof. \ep


\s{Proof of Theorem \ref{Th1}}
\renewcommand{\theequation}{4.\arabic{equation}}

\noindent
Now for $\al,d>0$, define
$$
\G_\la^\al:=\{u\in H_r^1(\R^3): \G_\la(u)\le\al\}
$$
and
$$
\mathcal{S}^d=\left\{u\in H_r^1(\R^3):  \inf_{v\in
\mathcal{S}_r}\|u-v\|\le d\right\}.
$$
Obviously, $\mathcal{S}_r\subset \mathcal{S}^d$, i.e.,
$\mathcal{S}^d\not=\phi$ for all $d>0$. For some $0<d<1$, we will
find a solution $u\in\mathcal{S}^d$ of problem \re{q2} for
sufficiently small $\la>0$.
The following proposition is crucial to obtain a suitable
(PS)-sequence for $\G_\la$ and plays a key role in our proof.
Choose
\be\lab{dc}
0<d<\min\left\{\frac{1}{3}\left[\frac{3}{2}\S^3\kappa^{-1}\right]^{\frac{1}{4}},\sqrt{3b}\right\},
\ee
where $\kappa$ is given in \re{f}.

\bo\lab{bo5} Let $\{\la_i\}_{i=1}^\infty$ be such that $\lim_{i\rg
\infty}\la_i=0$ and for all $i$, $\{u_{\la_i}\}\subset \mathcal{S}^d$
with
$$
\lim_{i\rg\infty}\G_{\la_i}(u_{\la_i})\le b\ \mbox{and}\
\lim_{i\rg\infty}\G_{\la_i}^{'}(u_{\la_i})=0.
$$
Then for $d$ small enough, there is $u_0\in \mathcal{S}_r$, up to
a subsequence, such that $u_{\la_i}\rg u_0$ in $H_r^1(\R^3)$.
\eo

\bp
For convenience, we write $\la$ for $\la_i$. Since $u_\la\in
\mathcal{S}^d$, there exist $U_\la\in\mathcal{S}_r$ and
$v_\la\in H^1(\R^3)$ such that $u_\la=U_\la+v_\la$ with
$\|v_\la\|\le d$. Since $\mathcal{S}_r$ is compact, up to a
subsequence, there exist $U_0\in\mathcal{S}_r$ and $v_0\in
H^1(\R^3)$, such that $U_\la\rg U_0$ strongly in $H^1(\R^3)$,
$v_\la\rg v_0$ weakly in $H^1(\R^3)$, $\|v_0\|\le d$ and $v_\la\rg v_0$ a.e. in
$\R^3$. Let $u_0=U_0+v_0$, then $u_0\in \mathcal{S}^d$ and $u_\la\rg
u_0$ weakly in $H^1(\R^3)$. It follows from
$\lim_{i\rg\infty}\G_{\la}^{'}(u_{\la})=0$ that $I'(u_0)=0$.
Now, we show $u_0\not\equiv 0$. Otherwise, if $u_0\equiv0$, then $\|U_0\|=\|v_0\|\le d$.
By \re{dc}, $\|\na U_0\|<\sqrt{3b}$. On the other hand, by $U_0\in\mathcal{S}_r$ and the Pohozaev's identity, $\|\na U_0\|=\sqrt{3b}$, which is a contradiction.
So $u_0\not\equiv0$ and $I(u_0)\ge b$. Meanwhile, thanks to
Lemma \ref{l2}, $\G_\la(u_\la)=I(u_0)+I(u_\la-u_0)+o(1)$, then we
have $I(u_\la-u_0)\le o(1)$. Thus, by \re{f} and the Sobolev'
embedding theorem, $\|u_\la-u_0\|^2\le\frac{2}{3}\kappa
\S^{-3}\|u_\la-u_0\|^6+o(1)$. If $\|u_\la-u_0\|\not\rg 0$ as $\la\rg
0$, up to a subsequence, we can get that $\|u_\la-u_0\|\ge
\left[\frac{3}{2}\S^3\kappa^{-1}\right]^{\frac{1}{4}}$ for $\la$
small. This is a contradiction. Thus, $u_\la\rg u_0$ strongly in
$H_r^1(\R^3)$. The proof is completed. \ep

\noindent
By Proposition \ref{bo5}, there exist
$$0<d<\min\left\{1,\frac{1}{3}\left(\frac{3}{2}\S^3\kappa^{-1}\right)^{\frac{1}{4}},\sqrt{3b}\right\}$$ and
$\omega>0, \la_0>0$ such that $\|\G_\la'(u)\|\ge\omega$ for
$u\in\G_\la^{D_\la}\bigcap(\mathcal{S}^d\setminus\mathcal{S}^{\frac{d}{2}})$
and $\la\in(0,\la_0)$. Then, we have the following proposition.

\bo\lab{bo6} There exists $\al>0$ such that for small $\la>0$,
$$
\G_\la(\g(s))\ge C_\la-\al\ \ \mbox{implies that}\ \ \g(s)\in
\mathcal{S}^{\frac{d}{2}},
$$
where $\g(s)=U(\frac{\cdot}{s}), s\in (0,t_0]$. \eo

\bp From a change of variables and the Pohoz\v{a}ev identity,
$$
\G_\la(\g(s))=\left(\frac{s}{2}-\frac{s^3}{6}\right)\int_{\R^3}|\na
U|^2\, \ud x+\frac{\la}{4}\int_{\R^3}\phi_{\g(s)}|\g(s)|^2\, \ud x.
$$
It follows from Lemma \ref{l0} that
$\G_\la(\g(s))=\big(\frac{s}{2}-\frac{s^3}{6}\big)\int_{\R^3}|\na
U|^2\, \ud x+O(\la^2)$. Note that
$$
\max_{s\in
[0,t_0]}\left(\frac{s}{2}-\frac{s^3}{6}\right)\int_{\R^3}|\na U|^2\, \ud x=b
$$
and $C_\la\rg b$ as $\la\rg 0$, the conclusion follows. \ep

\noindent
The next proposition assures the existence of a bounded
Palais-Smale sequence for $\G_\la$. Choosing small $\alpha_0>0$ satisfies
\be\lab{alpha}
\alpha_0\le\min\{\frac{\alpha}{2},\frac{1}{9}d\omega^2\}.
\ee
Noting that $\lim_{\la\rg 0}C_\la=\lim_{\la\rg 0}D_\la=b$, without
loss of generality, we can assume that $D_\la<C_\la+\alpha_0\le 2C_\la$
for $\la>0$ small enough.

\bo\lab{bo7} For any $\la>0$ small enough, there exists
$\{u_n\}_n\subset \G_\la^{D_\la}\cap \mathcal{S}^d$ such that
$\G_\la'(u_n)\rg 0$ as $n\rg\iy$.
\eo

\bp Assume by contradiction, there exists $a(\la)>0$ such that
$|\G_\la'(u)|\ge a(\la),\ u\in \mathcal{S}^d\cap \G_\la^{D_\la}$ for
some small $\la>0$. Then there exists a pseudo-gradient vector field
$T_\la$ in $H_r^1(\R^3)$ on a neighborhood $Z_\la$ of
$\mathcal{S}^d\cap \G_\la^{D_\la}$ (cf.\cite{Struwe}), such that for
any $u\in Z_\la$ satisfies
\begin{align*}
&\|T_\la(u)\|\le 2\min\{1,|\G_\la'(u)|\},\\
&\left<\G_\la'(u),T_\la(u)\right>\ge\min\{1,|\G_\la'(u)|\}|\G_\la'(u)|.
\end{align*}
Let $\eta_\la$ be a Lipschiz continuous function on $H_r^1(\R^3)$
such that $0\le \eta_\la\le 1,\ \eta_\la(u)\equiv 1$ if $u\in
\mathcal{S}^d\cap \G_\la^{D_\la}$ and $\eta_\la(u)=0$ if $u\in
H_r^1(\R^3)\setminus Z_\la$. Let $\xi_\la$ be a Lipschiz continuous
function on $\R$ such that $0\le \xi_\la\le 1,\ \xi_\la(t)\equiv1$
for $|t-C_\la|\le\frac{\alpha}{2}$ and $\xi_\la(t)=0$ for
$|t-C_\la|\ge \alpha$. Let
\begin{eqnarray*}
e_\la(u)= \left\{
\begin{array}{ll}
-\eta_\la(u)\xi_\la(\G_\la(u))T_\la(u),\ \mbox{if}\ u\in Z_\la,\\
0,\ \ \ \ \ \ \ \ \ \ \ \ \ \ \ \ \ \ \ \ \ \ \ \ \ \ \ \ \
\mbox{if}\ u\in H_r^1(\R^3)\setminus Z_\la,
\end{array}
\right.
\end{eqnarray*}
then for any $u\in H_r^1(\R^3)$, the following initial value problem
\begin{eqnarray*}
\left\{
\begin{array}{ll}
\frac{\ud}{\ud t}\Phi_\la(u,t)=e_\la(\Phi_\la(u,t)),\\
\Phi_\la(u,0)=u
\end{array}
\right.
\end{eqnarray*}
exists a unique global solution $\Phi_\la:H_r^1(\R^3)\times
[0,\iy)\rg H_r^1(\R^3)$ which satisfies
\begin{itemize}
\item [(1)] $\Phi_\la(u,t)=u$, if $t=0$ or $u\not\in Z_\la$ or $|\G_\la(u)-C_\la|\ge\alpha$,
\item [(2)] $\|\frac{\ud}{\ud t}\Phi_\la(u,t)\|_\la\le 2$, for all $u,t$,
\item [(3)] $\frac{\ud}{\ud t}\G_\la(\Phi_\la(u,t))\le 0$, for all $u,t$.
\end{itemize}

\noindent
With arguments similar as those in \cite{Chen-Zou}, we for any $s\in [0,t_0]$, there
exists $t_s\ge 0$ such that
$$
\Phi_\la(\g(s),t_s)\in \G_\la^{C_\la-\alpha_0},
$$
where $\g$ is given in Proposition \ref{bo6} and $\alpha_0$ is given in \re{alpha}. Let
$$
T_1(s):=\inf\{t\ge 0: \Phi_\la(\g(s),t)\in
\G_\la^{C_\la-\alpha_0}\}
$$
and
$\g_0(s)=\Phi_\la(\g(s),T_1(s))$, then
$\g_0$ is well defined in $[0,t_0]$ and there holds $\G_\la(\g_0(s))\le C_\la-\alpha_0$ for $s\in [0,t_0]$. With the similar arguments in
\cite{Byeon-Zhang-Zou,Chen-Zou}, we can get that $\g_0(s)$ is
continuous in $[0,t_0]$ and $\|\g_0(s)\|\le
\mathcal{C}_0+d<\mathcal{C}_0+1$. Thus, $\g_0\in \Upsilon_\la$ with
$\max_{t\in [0,t_0]}\G_\la(\g_0(t))\le C_\la-\alpha_0$, which is in
contradiction with the definition of $C_\la$. Therefore, the proof
is completed. \ep

\vskip2pt
\noindent

\noindent{\bf Proof of Theorem \ref{Th1} concluded.} It follows from Proposition \ref{bo7} that there exists
$\la_0>0$ such that for $\la\in(0,\la_0)$, there exists $\{u_n\}\in
\G_\la^{D_\la}\cap \mathcal{S}^d$ with $\G_\la'(u_n)\rg 0$ as
$n\rg\iy$. Assume that $u_n\rg u_\la$ weakly in $H^1(\R^3)$, then
$\G_\la'(u_\la)=0$. By the compactness of $\mathcal{S}_r$, we get that
$u_\la\in\mathcal{S}^d$ and $\|u_n-u_\la\|\le 3d$ for $n$ large.
In light of Lemma \ref{l0} and Lemma \ref{l2},
$$
\G_\la(u_n)=\G_\la(u_\la)+\G_\la(u_n-u_\la)+o(1).
$$
By the choice of $d$, it is easy to verify that
$\G_\la(u_n-u_\la)\ge0$ for large $n$. So, $\G_\la(u_\la)\le D_\la$.
Then $u_\la\in
\G_\la^{D_\la}\cap \mathcal{S}^d$ with $\G_\la'(u_\la)=0$. On the other hand, it is easy to know that $0\not\in\mathcal{S}^d$ for small $d>0$. Thus choosing $d>0$ small enough, $u_\la$ is a nontrivial solution of \re{q1}. In the following, we consider the asymptotic behavior of $u_\la$ as
$\la\rg 0$. Observe that
$$
\G_\la(u_\la)=I(u_\la)+\frac{\la}{4}\int_{\R^3}\phi_{u_\la}u_\la^2\,
\ud x,
$$
and that, for any $\varphi\in C_0^\iy(\R^3)$,
\begin{align*}
\G_\la'(u_\la)\varphi=I'(u_\la)\varphi+\la\int_{\R^3}\phi_{u_\la}u_\la\varphi\, \ud x.
\end{align*}
Note that $u_\la\in \mathcal{S}^d$.
Then, by Lemma \ref{l0} we get that
\begin{align*}
\hbox{$I(u_\la)\le D_\la$ and $I'(u_\la)\rg 0$ as $\la\rg 0$.}
\end{align*}
Assume that $u_\la\rg u$ weakly in $H_r^1(\R^3)$, then $I'(u)=0$.
Recall that $D_\la\rg b$ as $\la\rg 0$ and
$b\in(0,\frac{1}{3}\S^{\frac{3}{2}})$. It is easy to verify that
$u_\la\rg u$ strongly in $H_r^1(\R^3)$ as $\la\rg 0$ and $I(u)\le
b$. On the other hand, it follows from $u_\la\in \mathcal{S}^d$ that
$u\in \mathcal{S}^d$. Obviously, $0\not\in \mathcal{S}^d$ for $d$
small enough. Hence, choosing $d>0$ small enough, $u\not\equiv0$ and $I(u)\ge b$. Therefore,
$I(u)=b$, namely, $u$ is a least energy solution of problem~\re{q3}. \qed

\bigskip

\noindent{\bf Acknowledgements.}\! The authors would like to express their sincere gratitude to the anonymous referee for his/her careful reading and valuable suggestions. The first author also would like to thank Professor Zhi-Qiang Wang for his kind support and for fruitful discussions.

\end{document}